\documentclass[12pt]{article}
\usepackage{bbm, fancyhdr}
\usepackage{amssymb}
\usepackage{mathrsfs}
\usepackage{mathtools}
\usepackage{amsmath}
\usepackage{amsthm}
\usepackage{todonotes}
\usepackage[utf8]{inputenc}
\usepackage[english]{babel}
\usepackage{comment}
\usepackage{listings}
\usepackage{enumitem}
\usepackage{scalerel}
\usepackage{stackengine,wasysym}
\usepackage[backref=page]{hyperref}
\renewcommand*{\backref}[1]{}
\renewcommand*{\backrefalt}[4]{%
    \ifcase #1 (Not cited.)%
    \or        (Cited on page~#2.)%
    \else      (Cited on pages~#2.)%
    \fi}
\usepackage[utf8]{inputenc}
\usepackage{float}

\renewcommand{\leftmark}{{\small Complex structures on the product of two Sasakian manifolds}}

\hypersetup{colorlinks,
            breaklinks,
            urlcolor=magenta,
            linkcolor=magenta}

\usepackage{amsthm}
\usepackage{cleveref}

\newcommand\R{\mathbb{R}}
\newcommand\C{\mathbb{C}}

\newcommand\6{\partial}

\newcommand\End{\textrm{End}}

\newcommand\Id{\textrm{Id}}
\newcommand{\bas}{\textrm{bas}}
\newcommand{\Sym}{\textrm{Sym}}
\newcommand{\inv}{\textrm{inv}}
\renewcommand\Im{\rm{Im}}
\newcommand\im{\rm{im}}
\renewcommand\Re{\textrm{Re}}
\newcommand\Iso{\textrm{Iso}}
\newcommand\Lie{\textrm{Lie}}
\newcommand\pr{\textrm{pr}}

\newcommand\cali[1]{\mathcal{#1}}

\newcounter{globcounter}[section]

\newtheorem{theorem}[globcounter]{Theorem}
\newtheorem{proposition}[globcounter]{Proposition}

\newtheorem{lemma}[globcounter]{Lemma}
\theoremstyle{definition}
\newtheorem{remark}[globcounter]{Remark}
\theoremstyle{definition}
\newtheorem{claim}[globcounter]{Claim}
\theoremstyle{definition}
\newtheorem{definition}[globcounter]{Definition}

\newcounter{proofstepcounter}[globcounter]
\begin{document}

	\begin{center}
	{\Large\bf  Complex structures on the product of two Sasakian manifolds}\\[5mm]
	{\large Vlad Marchidanu\footnote{Partially supported by Romanian Ministry of Education and Research, Program PN-III, Project number PN-III-P4-ID-PCE-2020-0025, Contract  30/04.02.2021\\[1mm]
		
		
		}}\\[4mm]
	\end{center}
	
	{\small
		\hspace{0.15\linewidth}
		\begin{minipage}[t]{0.7\linewidth}
			{\bf Abstract.} A Sasakian manifold is a Riemannian manifold whose metric cone admits a certain K\"ahler structure which behaves well under homotheties. We show that the product of two compact Sasakian manifolds admits a family of complex structures indexed by a complex nonreal parameter, none of whose members admits any compatible locally conformally K\"ahler metrics if both Sasakian manifolds are of dimension greater than $1$. We compare this family with another family of complex structures which has been studied in the literature. We compute the Dolbeault cohomology groups of these products of compact Sasakian manifolds.
   \\
		\end{minipage}
	}

\tableofcontents
 
\section{Introduction}

Sasakian manifolds are the natural odd-dimensional analogue of K\"ahler manifolds (see e.g. \cite{boyer}). 
In the compact case, they are closely related to both projective and Vaisman manifolds (\cite{lck_book}).

K\"ahler manifolds can be viewed as  almost complex manifolds endowed with a Hermitian metric such that the associated fundamental two-form is parallel with respect to the metric connection.

Likewise, Sasakian manifolds can be thought of as almost contact manifolds endowed with a compatible Riemannian metric satisfying certain tensorial conditions (see \cite{blair} and Section \ref{subsubsection tensorial definition sasakians}).

Being even dimensional, a product of Sasakian manifolds is susceptible to bear almost complex structures. Indeed, more generally, Morimoto constructed an almost complex structure on the product of two almost contact manifolds (\cite{mori}) which proved to be integrable when the the two almost contact structures were normal.  If one starts with metric almost contact structures, then the product metric is compatible with Morimoto's almost complex structure. One thus obtains an almost Hermitian structure on the product. The usual complex structure of the Calabi-Eckmann manifold can be viewed this way. In particular, starting with two Sasakian manifolds (whose subjacent almost contact structures are normal), one obtains a Hermitian metric on the product.  

This product Hermitian structure on a Calabi-Eckmann manifold was later included by Tsukada in a two parameter family of Hermitian  structures (\cite{tsuk}). This construction was further generalized in \cite{lee} to the product of two Sasakian manifolds. It was recently considered also in \cite{andr_tolc}.

All these constructions use the tensorial definition of Sasakian manifolds and are heavily computational. With these techniques, the authors of \cite{andr_tolc} can prove that the considered two-parameter family of Hermitian structures is neither K\"ahler nor locally conformally K\"ahler.

What we propose in the present paper is a shift towards the modern definition of a Sasakian manifold as a Riemannian manifold with a K\"ahler structure on its Riemannian cone (\ref{subsubsection definition of sasakian via cone}). On the product of two compact Sasakian manifolds we construct a natural family of complex structures indexed by a purely complex parameter which we can prove that does not support neither K\"ahler nor locally conformally K\"ahler metrics. Furthermore, we are able to characterize the complex submanifolds of the product. Moreover, we show that our family of complex structures does not coincide, in general, with the one in \cite{lee}. Furthermore, we compute the Dolbeault cohomology groups of these complex manifolds.

{\bf Acknowledgements.} This paper is largely a result of my stay in Rio de Janeiro. I have learned many things from mathematicians at IMPA, but I would like to thank in the first place Prof. Misha Verbitsky for the insightful discussions we had and for helping me better understand mathematics. I am grateful to my thesis advisor, Prof. Liviu Ornea, who has guided me constantly and helped me with valuable comments and recommendations. 

\section{Sasakian manifolds}

\subsection{Tensorial definition of Sasakian manifolds}\label{subsubsection tensorial definition sasakians}
The notion of a Sasakian manifold was initially introduced by Shigeo Sasaki in \cite{sasaki} as an odd-dimensional counterpart to K\"ahler manifolds. We recall the tensorial definition of a Sasakian structure. 

Given a smooth, odd-dimensional manifold $S$, a Sasakian structure is given by the data $(g, \eta, \varphi, \xi)$, where $g$ is a Riemannian metric on $S$, $\eta$ is a $1$-form, $\varphi$ is a $(1,1)$-tensor field and $\xi$ is a vector field, satisfying the following properties for any $X, Y \in TS$:
\begin{equation*}
	\begin{split}
    \eta \circ \varphi &= 0 \\
    \eta(X) &= g(X,\xi)  \\
    \varphi^2 &= -\Id + \eta \otimes \xi \\
     g(\varphi(X),\varphi(Y)) &= g(X,Y) - \eta(X)\eta(Y) \\
     (-2d\eta \otimes \xi)(X,Y) &= \varphi^2([X,Y]) - \varphi([\varphi X, Y]) - \varphi([X,\varphi Y]) - [\varphi X, \varphi Y]
     \label{equation nijenhuis almost contact} \\
    \Lie_g \xi &= 0  \\
    (\nabla_X^g \varphi) Y &= g(X,Y) \xi - \eta(Y) X 
    \end{split}
\end{equation*}
A well-studied generalisation is that of an \emph{almost contact structure} $(S, \eta, \varphi, \xi)$, which occurs if we omit the presence of the metric $g$ and keep the first three conditions above, replacing $\eta(X) = g(X, \xi)$ with $\eta(X) = 1$, $\varphi(\xi) = 0$. This is usually viewed as a counterpart of almost complex geometry. See \cite{blair} for details.

In this paper we shall use the modern definition which places Sasakian manifolds into the framework of holonomy. This approach was widely spread following the pioneering work of C.P. Boyer and K. Galicki (see \cite{boyer}).

\subsection{Sasakian manifolds via the Riemannian cone} \label{subsubsection definition of sasakian via cone}
Let $(S,g)$ be an odd-dimensional Riemannian manifold and $(C(S):=(S \times \R^{>0}, g_{C(S)}= dt \otimes dt + t^2 g)$, $t\in\R^{>0}$, its Riemannian cone. 

\begin{definition}
	A Sasakian structure is the data of a K\"ahler structure $(J, \omega, g_{C(S)})$ on $C(S)$ such that the homothety map $h_\lambda: C(S) \rightarrow C(S)$, $h_\lambda(p,t):= (p, \lambda t)$ is holomorphic and satisfies $h_\lambda^* \omega = \lambda^2 \omega$ for each $\lambda \in \R^{>0}$.
\end{definition}

We denote by $R := t \frac{d}{dt}$ the Euler field on $C(S)$ and by $\xi:= JR$ the Reeb field. By definition $R$ is holomorphic, so $[R, \xi] = 0$. Since $C(S)$ is K\"ahler, $\xi$ is also holomorphic. When referring to $S$, we also denote by $\xi$ the vector field $\xi\rvert_{t=1}$ on $S \times \{1 \} \subset C(S)$.

The equivalence of the definition of Sasakian manifolds via their metric cone with the definition formulated in Subsection \ref{subsubsection tensorial definition sasakians} is established in \cite[Section~6.5]{boyer}. For our purposes, we mention that starting with a Sasakian manifold in the above sense, one defines the tensor field $\varphi \in \End(TS)$:
\[ \varphi(X):= \pr_{TS} JX, \quad X \in TS \subset TC(S) \]
where $J$ is the complex structure on $C(S)$.

We also define the $1$-form on $C(S)$, $\eta:= \frac{1}{t} J dt$, which is readily seen to be equal to $\frac{1}{t^2} i_R \omega$. As we did with $\xi$, we shall also denote $\eta$ the restriction $\eta = \eta\rvert_{t=1}$ on $S$. 
Then we have:

\begin{proposition} \label{claim square of varphi tensor}
$S$ is an almost contact manifold with contact form $\eta$ and characteristic field $\xi$. Moreover, ${\varphi^2 = -\Id + \eta \otimes \xi}$.
\end{proposition}

Denote by $D = \langle R, \xi \rangle^\perp$ the distribution ${g_{C(S)}}$-orthogonal to $\langle R, \xi \rangle$ on $C(S)$. Note that $t^2$ is a K\"ahler potential for $\omega$ and $dd^c (\log t)$ vanishes on $\langle R, \xi \rangle$, the  rest of its eigenvalues being positive. It follows that:
\begin{proposition} \label{proposition deta is kaehler on orthogonal}
   $\ker (d\eta) = \langle R, \xi \rangle$ and  $(d\eta)\rvert_D = \omega\rvert_{D} $. In particular $(d\eta)\rvert_D$ is a K\"ahler form.
\end{proposition}

\subsection{Basic cohomology of Sasakian manifolds}

\begin{definition}{\cite[Chapter~4]{tondeur}}\label{definition basic forms}
    Let $(M, \cali F)$ be a foliated manifold and consider $F \subset TM$ to be the subbundle of vectors tangent to leaves of $\cali F$. A form $\eta \in \Lambda^* M$ is called \textbf{basic} (with respect to $\cali F$) if for any vector field $X \in \Gamma F$, $\Lie_X \eta = 0$ and $i_X \eta = 0$.
\end{definition}

Denote the space of basic forms on a foliated manifold $(M,\cali F, F)$ by $\Lambda^*_\bas M$. By Cartan's formula, the exterior differential $d$ maps basic forms to basic forms. Therefore, $d$ induces a cohomology on basic forms, which we denote $H_\bas^* M$.

We are interested in a particular type of foliations:

\begin{definition}\label{definition transversally kaehler foliation}
        Let $(M,\cali F, F)$ be a foliated manifold.

        Let $\omega_0 \in \Lambda_{\bas}^*(M)$ with $d\omega_0 = 0$ and $g_0 \in \Sym_{\bas}^2(T^*M)$ such that $\omega_0|_F = 0$ and $\omega_0, g_0$ are positive definite on $TM/F$. 

        If the complex structure $J$ obtained from $\omega_0$ and $g_0$ is locally integrable on any open set in $M$ where the leaf space is defined, $(M,F, g_0, \omega_0)$ is called a \textbf{transversally Kähler foliation}.
\end{definition}

On compact K\"ahler manifolds the following well-known consequence of Hodge decomposition and Dolbeault decomposition holds.

\begin{theorem} {\rm(\cite[Theorem~VI.8.5]{demailly_complex_anal})} \label{theorem hodge decomposition}
Let $M$ be a compact K\"ahler manifold.
Denote by $H^{p,q}_{\bar \6} M$ the Dolbeault cohomology groups given by \newline
${\bar\6: \Lambda^{p,q} M \rightarrow \Lambda^{p,q+1} M}$. Then the Hodge decomposition holds:
    \[ H^k_{DR} M = \bigoplus\limits_{p+q = k} H^{p,q}_{\bar\6} M  \]
\end{theorem}

The usefulness of transversally K\"ahler foliations lies in the following result analogous to Theorem \ref{theorem hodge decomposition}.

\begin{theorem} {\rm(\cite[Theorem~30.28]{lck_book})} \label{theorem trasnversally kaehler admits hodge}
        Let $M$ be a compact manifold with a transversally Kähler foliation $(M,F, g_0, \omega_0)$ such that $F$ is generated by Killing vector fields and $M$ is equipped with a metric $g$ with $g|_{TM/F} = g_0|_{TM/F}$. Suppose  there exists $\Phi \in \Lambda^*(M)$ with $d\Phi = 0$, $\Phi|_F = 0$ and $\Phi$ is a volume form on $TM/F$.

        Then $H_{\bas}(M)$ behaves just like the cohomology of a K\"ahler manifold with respect to the K\"ahler form $\omega_0$. In particular, $H_{\bas}(M)$ admits the Hodge decomposition i.e.
        \[ H^k_{\bas} M = \bigoplus\limits_{p+q = k} H^{p,q}_{\bar\6_\bas} M  \]
        where $\bar\6_\bas$ is the operator given locally on the leaves of the foliation $F$ by the complex structure $J$ determined by $\omega_0$ and $g_0$ as in Definition \ref{definition transversally kaehler foliation}.
\end{theorem}

It turns out moreover that the cohomology of Sasakian manifolds is closely related to the basic cohomology of their associated transversally Kähler foliation. More precisely, we have:

    \begin{theorem} {\rm (\cite[Proposition~7.4.13]{boyer})}\label{theorem link between sasakian cohom and trasnversal cohom}
        Let $S$ be a Sasakian manifold of dimension $2n+1$ with characteristic (Reeb) field $\xi$. Let $F = \langle \xi \rangle$ be the transversally Kähler foliation generated by the Reeb field, which satisfies the conditions of the previous theorem. Then:
        \[ H^k(S) = \frac{H^k_{\bas}(S)}{\Im (\omega_0 \wedge \cdot)},\quad k < n \]
    \end{theorem}

\subsection{The product of two Sasakian manifolds}
In the context of (almost) contact geometry, Morimoto was the first to introduce an almost complex structure on the product of two almost contact manifolds in \cite{mori}. He shows that this almost complex structure is integrable if and only if condition \eqref{equation nijenhuis almost contact} is satisfied for each factor of the product. Building on Morimoto's ideas, Tsukada introduced in \cite{tsuk} a family of complex structures indexed by a complex nonreal parameter on the product of odd-dimensional spheres, noting that by the same argument as in \cite{mori} these structures are all integrable. In the same paper, Hermitian metrics associated with each of these complex structures are introduced. Watson generalised this family of pairs of complex structures and Hermitian metrics to products of Sasakian manifolds (\cite{watson}). We recall the definition of this family below. In the nomenclature of \cite{watson}, we call a structure in this family a Calabi-Eckmann-Morimoto structure, or CEM for short.

Let $S_1$, $S_2$ be Sasakian manifolds with $(1,1)$ tensors $\varphi_1$, $\varphi_2$ and Reeb fields $\xi_1, \xi_2$ respectively. Then there is a family of complex structures $\{ J_{a,b} : a, b \in \R, b \neq 0 \}$ on $S_1 \times S_2$:
\begin{align} 
J_{a,b}(X_1 + X_2)&:= \varphi_1(X_1) - \left( \frac{a}{b} \eta_1(X_1) + \frac{a^2+b^2}{b} \eta_2(X_2) \right) \xi_1 \nonumber \\
&+ \varphi_2(X_2) + \left( \frac{1}{b} \eta_1(X_1) + \frac{a}{b} \eta_2(X_2) \right) \xi_2
\label{calabi eckmann morimoto complex structure definition}
\end{align}

Let $g_i$ denote the Riemannian metric on the Sasakian manifold $S_i$, $i = \overline{1,2}$. For each pair $(a,b), b \neq 0$ there is an associated Hermitian (\cite{tsuk}) metric $g_{a,b}$ given by
\begin{align} 
g_{a,b}(X_1 + X_2, Y_1 + Y_2)&:= g_1(X_1, Y_1) + a \eta_1(X_1)\eta_2(Y_2) + a \eta_1(Y_1)\eta_2(X_2) \nonumber \\
&+ (a^2+b^2-1)\eta_2(X_2)\eta_2(Y_2)  + g_2(X_2, Y_2)
\label{calabi eckmann morimoto hermitian metric definition}
\end{align}

The metric data given by $g_{a,b}$ has been studied. It is shown in \cite[Theorem~1]{lee} that the metric $g_{a,b}$ is Einstein if and only if $a = 0$, $S_1$ is Einstein, and $S_2$ is $\eta_2$-Einstein with some specific constants (see \cite{okumura}, \cite{sasaki} for $\eta$-Einstein manifolds).
The authors of \cite{lee} also consider the property of weak $*$-Einsteiniainty for the product, which involves the interplay of $J_{a,b}$ with $g_{a,b}$. In showing that $(S_1 \times S_2, J_{a,b}, g_{a,b})$ is never weakly $*$-Einstein, they also prove that $(J_{a,b}, g_{a,b})$ is never K\"ahler.

Further exploring this interplay, the authors of \cite{andr_tolc} study whether and when the pairs $(J_{a,b}, g_{a,b})$ satisfy a number of natural conditions which are weaker than K\"ahlerianity, building on previous work in \cite{matsuo} and \cite{fino_ugarte}. The results known about $(J_{a,b}, g_{a,b})$ are summarized in the following

\begin{theorem}{\rm(\cite{matsuo},\cite{fino_ugarte},\cite{andr_tolc})}
    Let $S_1$ and $S_2$ be Sasakian manifolds of dimensions $2n_1+1$ and $2n_2+1$ respectively. Consider the complex structure $J_{a,b}$ (\ref{calabi eckmann morimoto complex structure definition}) and the Hermitian metric $g_{a,b}$ (\ref{calabi eckmann morimoto hermitian metric definition}). Then:
    \begin{enumerate}
        \item If $n_1 + n_2 \geq 1$ then $(J_{a,b}, g_{a,b})$ is not balanced (see \cite{michelson balanced}).
        \item  $(J_{a,b}, g_{a,b})$ is LCK (see \cite{vaisman_lck} as well as \cite[Chapter~3]{lck_book} for equivalent definitions) if and only if $n_1 = 0$ and $n_2 \geq 1$ or $n_2 = 0$ and $n_2 \geq 1$; if it is LCK, then it is also Vaisman.
        \item  $(J_{a,b}, g_{a,b})$ is SKT (see \cite{howe_papadopoulos_skt}) if and only if either $n_1 = 1$ and $n_2 = 0$ or $n_1 = 0$ and $n_2=1$ or $a = 0$ and $n_1 = n_2 = 1$.
        \item If $n_1+n_2 \geq 2$ then the condition
        \begin{equation*}
             n_1(n_1-1) + 2an_1 n_2 + n_2(n_2-1)(a^2+b^2) = 0
        \end{equation*}
        holds if and only if $(J_{a,b}, g_{a,b})$ is $1$-Gauduchon (see \cite{fu_k_gauduchon} for $k$-Gauduchon) if and only if $(J_{a,b}, g_{a,b})$ is astheno-K\"ahler (see \cite{jost_yau_astheno_kaehler}). 
        \item If $n_1 + n_2 \geq 3$ and $2 \leq k \leq \dim_\C (S_1 \times S_2) -1$, then $(J_{a,b}, g_{a,b})$ is $k$-Gauduchon if and only if the following holds:
        \begin{equation*}
            (n_1 + n_2 - k)\left( n_1(n_1-1) + 2an_1 n_2 + n_2(n_2-1)(a^2+b^2) \right) = 0
        \end{equation*}
    \end{enumerate}
\end{theorem}

\section{The main result}

\begin{theorem}\label{theorem main}
    Let $S_1$, $S_2$ be compact Sasakian manifolds of respective dimensions $2n_i + 1$, with $n_i > 1$. Then $S_1 \times S_2$ has a family of complex structures indexed by a complex nonreal parameter, none of whose members admits any Kähler or LCK metrics.
    
\begin{proof}
    \textbf{Step 1. Definition of the complex structure on the product.} 
    
    To define the complex structure, we consider the following generalisation of Calabi-Eckmann manifolds.
    Let $S$ be a Sasakian manifold and define an action of $(\C, +)$ on the open cone $C(S)$ by putting for $a+bi \in \C$:
    \begin{equation} \label{action_def}
    	(a+b\sqrt{-1})\cdot p := \phi^{aR + bJR}_1(p), 
    \end{equation}
 where $\phi^X_t$ denotes the flow of  the vector field $X$ at time $t$.
    Since $R$ and $\xi$ commute, we have 
    \[ \phi^{(cR+d\xi)+(aR+b\xi)}_1(p) = \phi^{cR+d\xi}_1(\phi^{aR+b\xi}_1(p)) \]
    In other words
    \[ (( c+d\sqrt{-1} )+ (a+b\sqrt{-1}))\cdot p = (c+\sqrt{-1}d)\cdot ( (a+b\sqrt{-1})\cdot p )\]
    showing that indeed \eqref{action_def}  defines a group action.

    This action is a holomorphic map $\C \times C(S) \rightarrow C(S)$. Indeed, the Reeb and Euler fields act by biholomorphisms. Further, let $x \in C(S)$ and $v \in \C$, $X \in T_v \C$ and $\gamma(t)$ be a curve with tangent vector $X$ at $v$. Then $JX = \frac{d}{dt}\rvert_{t=0} \left( \sqrt{-1} \gamma(t) \right)$. Since $[R, \xi] = 0$, one vector field is invariated by the flow of the other. Therefore:
   \begin{align*}
   d_v( w \mapsto (w \cdot x))\left( JX \right) &=
   \frac{d}{dt}\rvert_{t=0} \left( (\sqrt{-1} \gamma(t)) \cdot x \right) \\
    &=\frac{d}{dt}\rvert_{t=0} \left( \phi_1^{ (-\Im(\gamma(t)))R + \Re(\gamma(t)) \xi }(x) \right) \\
    &= \frac{d}{dt}\rvert_{t=0} \left( 
        \phi_{\Re(\gamma(t))}^{\xi}\left( \phi_{-\Im(\gamma(t))}^{R}(x) \right)
    \right) \\
    &= -\Im(X) R + \Re(X) \xi = J(\Re(X) R + \Im(X) \xi) \\
    &= J d_v( w \mapsto (w \cdot x))\left( X \right)
   \end{align*}
   which shows that the map $w \mapsto (w \cdot x)$ is holomorphic for every fixed $x \in C(S)$. 

    Let now $S_i$, $i=1,2$, be compact Sasakian manifolds with Euler fields $R_i$, Reeb fields $\xi_i:= J_i R_i$, and consider the diagonal action of $\C \times \C$ on the product of the cones $C(S_1) \times C(S_2)$.

    Fix some $\alpha \in \C$ with $\Im \alpha \neq 0$ and define the subgroup ${G_\alpha := \{ (t,\alpha t): t \in \C \}}$ of $(\C \times \C, +)$.  Clearly, $G_\alpha$ is isomorphic with $\C$ and acts on $C(S_1) \times C(S_2)$.  We analyze $(C(S_1) \times C(S_2))/G_\alpha$.

    Let $r_i: C(S_i) \rightarrow \R^{>0}$ be  the projections on the radial directions.

    \begin{claim}\label{claim calabi eckmann action on product sasaki cones}
        For any $(a,b) \in {\R^{>0}}\times \R^{>0} $ and any $x=(\tilde p_1, \tilde p_2) \in C(S_1) \times C(S_2)$ there exists a unique $v \in \C \simeq G_\alpha$ such that 
        $r_1(v \cdot x) = a$ and $r_2(v \cdot x) = b$.
        \begin{proof}
        For $\tilde p_i=(p_i, t_i)$, $p_i \in S_i$, $t_i \in \R^{>0}$ we have:
        \begin{align*}
        v \cdot p_1 &= \phi_1^{\Re( v ) R_1 + \Im (v) \xi_1}(\tilde p_1) = \phi_1^{\Im (v) \xi_1}(\phi_1^{\Re (v) R_1}(\tilde p_1)) \\
        &= \phi_1^{\Im (v) \xi_1}((p_1, e^{\Re (v)} t_1))
        \end{align*}
        Since $\Im (v) \xi_1$ acts only on the level sets of the cone, when we project to the radial direction we get:
        \[ r_1(v \cdot \tilde p_1) = e^{\Re (v)} t_1 \]
        Similarly
        $ r_2(\alpha v \cdot \tilde p_2) = e^{\Re (\alpha v)} t_2 $.
        Therefore, what we need to show is that for any strictly positive $a, b, t_1, t_2$ there exists a unique $v \in \C$ such that:
        \[ \begin{cases}
        e^{\Re (v)} t_1 = a \\
        e^{\Re (\alpha) \Re (v) - \Im (\alpha) \Im (v)} t_2 = b
        \end{cases} \]
        or
        \[
        \begin{cases}
        \Re (v) = \log(a) - \log(t_1)\\
        \Re (\alpha) \Re (v) - \Im (\alpha) \Im (v) = \log(b) - \log(t_2)
        \end{cases}
        \]
        Since $\Im \alpha \neq 0$, we have
        \begin{align*}
            \Im v = \frac{\Re (\alpha) (\log(a) - \log(t_1)) + \log(t_2) - \log(b)}{\Im (\alpha)},
        \end{align*}
        and hence the solution $v$ exists and is unique.
        
        \end{proof}
    \end{claim}

    Claim \ref{claim calabi eckmann action on product sasaki cones} provides an identification of the quotient $(C(S_1) \times C(S_2) )/G_\alpha$ with $S_1 \times S_2$, given by an explicit formula for $\pi: C(S_1) \times C(S_2) \longrightarrow S_1 \times S_2$. Denote, as in Claim \ref{claim calabi eckmann action on product sasaki cones},
    \begin{equation} \label{equation formula for v t_1 t_2}
    	v(t_1, t_2) := - \log(t_1) + \frac{\sqrt{-1}}{\Im (\alpha)} ( - \Re (\alpha) \log t_1 + \log t_2 ).
    \end{equation}
Now the map $\pi: C(S_1) \times C(S_2) \longrightarrow S_1 \times S_2$ can be described as:
    \begin{equation}\label{explicit formula for pi}
    \pi( (p_1, t_1), (p_2, t_2) ) = \left( \phi_1^{\Im (v(t_1,t_2)) \xi_1} (p_1),\  \phi_1^{\Im (\alpha v(t_1, t_2)) \xi_2}(p_2) \right).
    \end{equation}
    
    Since the action of $G_\alpha$ defines a holomorphic map 
    \[
    G_\alpha \times C(S_1) \times C(S_2) \rightarrow  C(S_1) \times C(S_2)
    \]
    
   We conclude that 
    \[ 
    (C(S_1) \times C(S_2))/G_\alpha \simeq S_1 \times S_2
    \]
    admits a complex structure compatible with the smooth product structure on $S_1 \times S_2$ and making the projection map $\pi: C(S_1) \times C(S_2) \longrightarrow S_1 \times S_2$ a holomorphic submersion.
    
\medskip

\textbf{Step 2.} We now aim to better understand the complex structure induced by $\pi$. More precisely, we show that on the transverse distributions of each Sasakian, it acts like the complex structure on the cone, while it takes each Reeb field to the span of the two Reeb fields.

To keep notation simple, we will deliberately use the same notation $\xi_i$ for the Reeb field(s) both on the product of the K\"ahler cones and on the product of the Sasakian manifolds.

Let $X \in T_{p_1} S_1$ and $x \in \pi^{-1}(p_1, p_2) \subset C(S_1) \times C(S_2)$ for some $p_2 \in S_2$. We see $X$ as tangent in $x$ to $C(S_1) \times C(S_2)$. 
$X$ can be extended to a vector field $\tilde X$ on $C(S_1) \times C(S_2)$, such that $\tilde X$ is tangent to $S_1$ and moreover $\tilde X$ commutes with $\xi_1$ (and hence with all multiples of $\xi_1$) in a neighborhood of $x$. We can obtain such an extension by considering a chart on $S_1$ in which $\xi_1$ is a standard coordinate vector field, extending the expression of $X$ in this chart to a constant vector field and multiplying it with a bump function.

An extension $\tilde X$ of $X$ with $[\tilde X, \xi_1 ] = 0$ guarantees that
${ d_x \phi_1^{\xi_1} (\tilde X_x) = \tilde X_{\phi_1^{\xi_1}(x)} }$

Further, we have:
\begin{align*}
    d_x \pi (\tilde X) &= \frac{d}{dt}\rvert_{t=0} \left( \pi( (\phi_t^{\tilde X}(p_1), t_1),\quad  (p_2, t_2) ) \right) \\
    &= \frac{d}{dt}\rvert_{t=0} \left( \phi_1^{\Im (v(t_1,t_2)) \xi_1} (\phi_t^{\tilde{X}}(p_1)), \quad \phi_1^{\Im (\alpha v(t_1,t_2)) \xi_2} (p_2) \right) \\
    &= \left( d_{p_1} \left( p \mapsto \phi_1^{\Im (v(t_1,t_2)) \xi_1}(p) \right)(\tilde{X}_{p_1} ), \quad 0 \right)_{\pi(x)} \\
    &=X
\end{align*}

Similarly, for $X \in T_{p_2} S_2$ we have:
\begin{align*}
    d_x \pi (\tilde X) 
    &= \frac{d}{dt}\rvert_{t=0} \left( \phi_1^{\Im v(t_1,t_2) \xi_1} (p_1), \quad \phi_1^{\Im (\alpha v(t_1, t_2)) \xi_2} (\phi_t^{\tilde{X}}(p_2)) \right) 
    = X
\end{align*}

For the first Euler field:
\begin{align*}
    d_{x} \pi( R_1 ) &=
    \frac{d}{dt}\rvert_{t=0} \pi \left( (p_1, e^t t_1), \quad (p_2, t_2) \right) \\
    &=  \frac{d}{dt}\rvert_{t=0} \left( \phi_1^{\Im (v(e^t t_1,t_2)) \xi_1} (p_1), \quad  \phi_1^{\Im (\alpha v(e^t t_1,t_2)) \xi_2} (p_2) \right) \\
    &=  \frac{d}{dt}\rvert_{t=0} \left( \phi_{\Im ( v(e^t t_1,t_2))}^{\xi_1} (p_1), \quad  \phi_{\Im (\alpha v(e^t t_1,t_2))}^{ \xi_2} (p_2) \right) \\
    &=  \frac{d}{dt}\rvert_{t=0} \left( \Im( v(e^t t_1,t_2)) \right) \left( \xi_1 \right)_{\pi(x)}
        +
        \frac{d}{dt}\rvert_{t=0} \left( \Im(\alpha v(e^t t_1,t_2)) \right) \left( \xi_2 \right)_{\pi(x)}
\end{align*}

Denote from now $a = \Re \alpha, b = \Im \alpha$.

According to \eqref{equation formula for v t_1 t_2}, $v(e^t t_1,t_2) = v(t_1,t_2) - t\left(1 + \frac{a}{b}\sqrt{-1} \right)$. Hence
    \begin{align*}
        d_{x} \pi( R_1 )
        &= - \frac{1}{b} \left( a    \xi_1  + (a^2+b^2) \xi_2 \right)_{\pi(x)}
    \end{align*}

For the second Euler field, since, by \eqref{equation formula for v t_1 t_2},
$v(t_1, e^t t_2) = v(t_1,t_2) + \frac{t}{b}\sqrt{-1}$,
we deduce as before:
\begin{align*}
    d_{x} \pi( R_2 ) &=
    \frac{d}{dt}\rvert_{t=0} \pi \left( (p_1, t_1), \quad (p_2, e^t t_2) \right) \\
    &=  \frac{d}{dt}\rvert_{t=0} \left( \phi_1^{\Im (v(t_1, e^t t_2)) \xi_1} (p_1), \quad  \phi_1^{\Im (\alpha v(t_1, e^t t_2)) \xi_2} (p_2) \right) \\
    &= \frac{1}{b} \left( \xi_1 + a \xi_2 \right)_{\pi(x)}
\end{align*}

In summary, we have:
\begin{align}
    d_x \pi( X ) &= X, \quad x=( (p_1, 1), (p_2, 1) ), \quad X \in T_{p_1} S_1 \sqcup T_{p_2} S_2
    \label{equation pi of vectors tangent to sasakians}
    \\
    d_x \pi( R_1 )
        &= - \frac{1}{b} \left( a   \xi_1 + (a^2+b^2)  \xi_2   \right)_{\pi(x)}
    ,\quad
    d_{x} \pi( R_2 ) = \frac{1}{b} \left( \xi_1 + a \xi_2 \right)_{\pi(x)}
    \label{equation pi of eulers}
\end{align}

\textbf{Step 3. The above family of complex structures does not admit any compatible Kähler metric}.

\textbf{Step 3.1.} Let $\eta_1$ be the pullback of the contact form on $S_1$ through $S_1 \times S_2 \rightarrow S_1$. Then $d\eta_1$ is a semipositive $(1,1)$-form.

Indeed, to see that $d\eta_1$ is $(1,1)$, it's enough to check that $d\eta_1(z \pi_* X, z \pi_* Y) = z\bar z d\eta_1(\pi_*X, \pi_*Y)$ for $z \in \C$. By \eqref{equation pi of vectors tangent to sasakians} and holomorphicity of $\pi$:
\begin{align*}
    d\eta_1(z \pi_* X, z \pi_* Y) =& \Re(z)^2 d\eta_1(  X,  Y) 
    + \Im(z)^2 d\eta_1(\pi_* JX, \pi_*JY) \\
    &+ \Re(z)\Im(z) \left( d\eta_1(X, \pi_* JY) + d\eta_1(\pi_* JX, Y) \right)
\end{align*}

Suppose $X$ is orthogonal to $\langle R, \xi \rangle$ on its respective Sasakian manifold. If $Y$ is also orthogonal, we are done since $\pi_* JY = JY$ and $d\eta_1$ is transversally K\"ahler on the cone. Otherwise $Y$ is a multiple of a Reeb vector $\xi_i$, so $\pi_* JY \in \langle \xi_1, \xi_2 \rangle$, so $\pi_* JY \in \ker d\eta_1$ and since also $Y \in \ker d\eta_1$, the wanted equality checks trivially. Finally, if $X$ is a multiple of a Reeb vector, the wanted equality checks trivially because again $\{ X, \pi_* JX \} \subset \ker d\eta_1$.

Now checking semipositivity is equivalent by the holomorphicity of $\pi$ to checking that for $X \in TC(S_1) \times TC(S_2)$ we have
$ d\eta_1( \pi_* X, \pi_* JX) \geq 0 $. 
If $X$ is tangent to either $S_1$ or $S_2$ and is transverse to the Euler and Reeb fields, then $JX$ stays outside the distribution generated by the Euler and Reeb fields, and so by \eqref{equation pi of vectors tangent to sasakians} $d\eta_1( \pi_* X, \pi_* JX) = d\eta_1(X, JX)$ and the latter is a nonegative quantity because $d\eta_1$ is semipositive on the cone.
If $X$ is either $\xi_1$ and $\xi_2$ then $d\eta_1( \pi_* X, \pi_* JX) = 0$ by \eqref{equation pi of eulers}.

\textbf{Step 3.2.}
        Suppose $S_1 \times S_2$ is Kähler with Kähler form $\omega$. 
        \[ d(\eta_1 \wedge \omega^{\dim_\C (S_1 \times S_2) - 1}) = (d\eta_1) \wedge \omega^{\dim_\C (S_1 \times S_2) - 1}  \]
        because $d\omega = 0$. So by Stokes' Theorem
        \begin{equation}\label{integral prod sasaki contact with Kähler} \int_{S_1 \times S_2} (d\eta_1) \wedge \omega^{\dim_\C (S_1 \times S_2) - 1} = 0
        \end{equation}
        Because $d \eta_1$ is semipositive, $d\eta_1 \wedge \omega^{\dim_\C (S_1 \times S_2) - 1}$ is a semipositive volume form, which vanishes if and only if $d\eta_1$ vanishes. But then $d\eta_1$ vanishes by \eqref{integral prod sasaki contact with Kähler}, which contradicts the fact that $d\eta_1$ is positive on the distribution transverse to $\ker d\eta_1$.
        
\textbf{Step 4.}
Let $S_1$, $S_2$ be  Sasakian manifolds of respective dimensions $2n_i + 1$ with $n_i > 1$. By the Künneth formula 
$H^1(S_1 \times S_2) = H^1(S_1) \oplus H^1(S_2)$.
In view of Theorem \ref{theorem link between sasakian cohom and trasnversal cohom}, we can represent forms in $H^1(S_i)$ with basic forms. Hence, in view of Theorem \ref{theorem trasnversally kaehler admits hodge}, we can represent $[\eta] \in H^1(S_1 \times S_2)$ as $\eta^{1,0} + \eta^{0,1}$, where $\eta^{1,0}$ is holomorphic and closed and $\eta^{0,1}$ is antiholomorphic and closed. To see that this is the case, suppose that $\alpha$ is a holomorphic representative of a basic class on one of the Sasakian manifolds, say $[\alpha] \in H^*_\bas (S_1)$. The fact that $\alpha$ is a basic holomorphic form implies that $\pi_1^* \alpha$ is holomorphic, where $\pi_1: C(S_1) \rightarrow S_1$ is the projection. We need to check that this implies that $\alpha$ is holomorphic as a form on $S_1 \times S_2$ with the complex structure induced by the projection $\pi$ from $C(S_1) \times C(S_2)$. But by \eqref{equation pi of vectors tangent to sasakians} and \eqref{equation pi of eulers}, we obtain $\pi^* \alpha = \pi_1^* \alpha$ up to a constant, and hence $\alpha$ is holomorphic.

\textbf{Step 5.}
            Assuming $S_1 \times S_2$ is LCK, we represent the Lee form $\theta$ as $\theta = \theta^{1,0} + \theta^{0,1}$ with $\theta^{1,0}$ holomorphic and closed and $\theta^{0,1}$ antiholomorphic and closed. Thus we get $dd^c \theta = 0$. Then $dd^c (\omega^{n-1}) = \omega^{n-1} \wedge \theta \wedge J\theta$, so 
            \[ \int_M \omega^{n-1} \wedge \theta \wedge J\theta = 0\]
            Combined with the fact that $\theta \wedge J\theta$ is semipositive $(1,1)$, the above equality shows that $\omega^{n-1} \wedge \theta \wedge J\theta = 0$. 

            So $\theta \wedge J\theta = 0$.
            Hence $\theta = 0$ since $\theta$ and $J\theta$ are linearly independent. This shows that $S_1 \times S_2$ is GCK, but then it also admits a Kähler structure, which is a contradiction by Step 3.
    \end{proof}
\end{theorem}

\begin{remark}
    \begin{sloppypar}
    The same proof as in Step 3 shows that $S_1 \times S_2$ does not admit balanced metrics i.e. metrics with Hermitian form $\omega$ satisfying $d \omega^{\dim_\C(S_1 \times S_2) -1 } = 0$, since in that case we also obtain that ${ (d\eta_1) \wedge \omega^{\dim_\C (S_1 \times S_2) - 1} }$ is exact.
    \end{sloppypar}
\end{remark}

\begin{remark}
    The argument developed in Steps 3 through 5 also shows that the CEM complex structure defined by \eqref{calabi eckmann morimoto complex structure definition} does not admit any compatible locally conformally K\"ahler metric.
\end{remark}

\section{Complex submanifolds of the product of Sasakian manifolds}

Let $S_1,S_2$ be compact Sasakian manifolds with $\dim_\R S_i = 2n_i+1$ and with contact forms $\eta_1, \eta_2$. Let $S_1 \times S_2$ be their product with the complex structure induced by the action of $G_\alpha$ on the product of their cones as in the proof of Step 1 of Theorem \ref{theorem main}. 

\begin{theorem}
	Let $Z \subset S_1 \times S_2$ be a complex submanifold of $\dim_\C Z = k$ where the complex structure on $S_1 \times S_2$ is induced by the Calabi-Eckmann action on the product of the cones. Then $Z$ is tangent to $\ker(d\eta_1 + d\eta_2)$.
	\begin{proof}
		Let $\eta = \eta_1 + \eta_2$. Then We have:
		\[ d(\eta \wedge (d\eta)^{k-1}) = d\eta \wedge (d\eta)^{k-1} = (d\eta)^k \]
		So by Stokes' theorem we have:
		\[ \int_Z (d\eta)^k = 0 \]
		Since outside $\ker(d\eta)$, $d\eta$ is strictly positive, and $Z$ is a complex submanifold, we must thus have $TZ \subset \ker(d\eta)$.
	\end{proof}
\end{theorem}

\section{Comparison with the \newline Calabi-Eckmann-Morimoto complex structures}

Consider again the principal $G_\alpha = \{ (v, \alpha v): v \in \C \}$-bundle 
\[ {\pi: C(S_1) \times C(S_2) \rightarrow S_1 \times S_2} \] 
where $\alpha \in \C \setminus \R$. 

The natural question arises whether $J_{a,b}$ defined by \eqref{calabi eckmann morimoto complex structure definition} coincides with the complex structure induced by $\pi$. 

\begin{theorem}
For every fixed $\alpha \in \C \setminus \R$, the complex structure induced by $G_\alpha$
 does not in general coincide with the complex structure $J_{a,b}$ for any $a, b \in \R$, $b \neq 0$.
\begin{proof}
 For general Sasakian manifolds $S_1$ and $S_2$, by uniqueness of the complex structure making ${\pi: C(S_1) \times C(S_2) \rightarrow S_1 \times S_2}$ a holomorphic submersion, the complex structures on $S_1 \times S_2$ coincide if and only if for any $x = (p_1,t_1, p_2,t_2)$ and any $X_i \in T_{(p_i,t_i)} C(S_i)$, $i = \overline{1,2}$, we have
 \begin{equation*} 
J_{a,b} d_{x} \pi (X_1 + X_2) = d_{x} \pi( J(X_1+X_2) )
\end{equation*}

For $X_1 = \xi_1, X_2 = 0$ we have by \eqref{equation pi of vectors tangent to sasakians}
\[ J_{a,b} \pi_* \xi_1 = J_{a,b} \xi_1 = \frac{1}{b}\left( -a \xi_1 + \xi_2 \right) \]
while by \eqref{equation pi of eulers}
\[ \pi_* (J\xi_1) = - \pi_*(R_1) = \frac{1}{\Im (\alpha)} \left( \Re (\alpha) \xi_1 + ( (\Re (\alpha))^2 + (\Im (\alpha))^2) \xi_2) \right) \]

Furthermore,
    \[ J_{a,b} \pi_* \xi_2 = \frac{1}{b} \left( -(a^2+b^2) \xi_1 + a \xi_2 \right) \]
and 
\[ \pi_* (J\xi_2) = {-\pi_* (R_2)}= {-\frac{1}{\Im (\alpha)} (\xi_1 + \Re (\alpha) \xi_2)}. \]
So if $J_{a,b}$ coincides with the structure induced by $\pi$ we obtain the following system of equations:

\[
\begin{cases}
\Re(\alpha) b = - a \Im(\alpha), \quad \Im(\alpha) = b (\Re (\alpha))^2 + (\Im (\alpha))^2) \\
b = (a^2 + b^2)\Im(\alpha), \quad \Re(\alpha)b = -a\Im(\alpha) 
\end{cases}
\]
This leads to the equation:
\[ (\Im (\alpha))^4 + (2(\Re (\alpha))^2 - 1) (\Im (\alpha))^2 + (\Re (\alpha))^4 = 0 \]
which implies that $|\Im(\alpha)| \leq 1$ and that $(\Re (\alpha))^2 \leq \frac{1}{4}$. Hence, for $\alpha$ such that these conditions are not met, we cannot find $(a,b)$ such that $J_{a,b}$ coincides with the complex structure induced by $\pi$.

However, on $\langle \xi_1, \xi_2 \rangle$ $J$ is $-J_{\Re \alpha, \Im\alpha}^T$ where the superscript is matrix transpose. 
\end{proof}
\end{theorem}

\section{The Dolbeault cohomology of the product of Sasakian manifolds}
Let $S_1, S_2$ be compact Sasakian manifolds with the action of $G_\alpha$ as in Step 1 of Theorem \ref{theorem main}, $\alpha = a+b\sqrt{-1}$, $a \in \R, b \in \R\setminus \{ 0 \}$. Denote from now $M:= S_1 \times S_2$.

 Consider $\eta_1, \eta_2$ the two contact forms on $M$. Let $\eta:= \eta_1 + \eta_2$, $\omega_0:= d\eta$ and $\eta^{0,1}, \eta^{1,0}$ be the $(0,1)$ and $(1,0)$ parts of $\eta$, respectively. Since $\omega_0$ is a $(1,1)$-form and $\omega_0 = d\eta = \6 \eta^{0,1} + \6 \eta^{1,0} + \bar\6 \eta^{0,1} + \bar\6 \eta^{1,0}$, we obtain that $\6 \eta^{1,0} = 0$ and $\bar \6 \eta^{0,1} =0$.

 Endow $M$ also with a Hermitian metric such that the two Reeb fields are Killing, as follows. Consider 
$V:=\langle \xi_1, \xi_2 \rangle$ with the frame $\{ \xi_1, \xi_2 \}$ and $J_{\alpha}$ the complex structure induced by $\pi_{\alpha}$ as in Step 1 of Theorem \ref{theorem main}. Recall also that $J_{a,b}$ is defined as in $\eqref{calabi eckmann morimoto complex structure definition}$ and the metric $g_{a,b}$ defined as in \eqref{calabi eckmann morimoto hermitian metric definition} is Hermitian with respect to $J_{a,b}$. On $V$ we have $J_{\alpha}|_V = -(J_{a,b}|_V)^T$ for $a=\Re(\alpha), b=\Im(\alpha)$. Hence $J_{\alpha}|_V$ is the negative of the morphism induced on $V^*$ by $J_{a,b}|_V$, so $J_{\alpha}|_V$ is Hermitian with respect to $(g_{a,b}|_V)^{-1}$. Since on $V^{\perp}$ $J_{a,b}$ coincides with $J_\alpha$, the metric

\begin{align*}
g_\alpha := & g_1+g_2 
 - ab^{-2} \left( \eta_1 \otimes \eta_2 +\eta_2 \otimes \eta_1 \right) \\
& +\left( b^{-2}(a^2+b^2)-1 \right) \eta_1 \otimes \eta_1
+ (b^{-2} -1) \eta_2 \otimes \eta_2
\end{align*}

is Hermitian on $M$ with respect to $J_{\alpha}$, where $\eta_i$ and $g_i$ are the contact forms and Riemannian metrics respectively on $S_i$, extended with $0$ on the Sasakian manifold they are not initially defined on. Moreover, $\xi_1, \xi_2$ are Killing with respect to $g_{\alpha}$ because each $\xi_i$ is Killing with respect to $g_i$ and $\Lie_{\xi_i} \eta_i = 0$ (because $S_i$ is contact with characteristic field $\xi_i$ and by Cartan's formula).
\begin{sloppypar}
  By a theorem of Myers and Steenrod (\cite{mayers_steenrod}), $\Iso_{g_{\alpha}}(M)$ is a Lie group, which is compact since both Sasakian manifolds are compact. Consider $K$ to be the closure of the subgroup generated by $\phi_t^{\xi_1}$ and $\phi_t^{\xi_2}$ inside $\Iso_{g_{\alpha}}(M)$. By the closed subgroup theorem, $K$ is also a (compact) Lie group. Take $\Lambda^*(M)^{\inv}$ to be all forms on $M$ which are invariant under $K$.
  A standard continuity argument shows that 
  ${(\Lambda^*(M))^{\inv} = \{ \alpha \in \Lambda^*(M): \Lie_{\xi_1} \alpha = \Lie_{\xi_2} \alpha = 0 \}}$.
  Consider also $(\Lambda^*(M))_\bas$ to be all the basic forms with respect to the foliation $\langle \xi_1, \xi_2 \rangle$; clearly $(\Lambda^*(M))_{\bas} \subset (\Lambda^*(M))^\inv$ (see Definition \ref{definition basic forms}). Locally, basic forms come from the leaf space of the foliation.
\end{sloppypar}

Put $\Lambda_{B, \eta^{0,1}}^{p,q}:= (\Lambda^{p,q})_\bas \oplus \left( \eta^{0,1} \wedge \Lambda^{p,q-1}_\bas  \right)$. 

Since $\bar\6 \eta^{0,1} = 0$, for each $p \geq 0$ the restriction of $\bar\6$ gives a complex 
    \[ \bar\6: \Lambda_{B, \eta^{0,1}}^{p,*} \rightarrow \Lambda_{B, \eta^{0,1}}^{p,*+1} .\]

Now consider the operator  $L_{\omega_0}: \Lambda^*(M) \rightarrow \Lambda^*(M)$ to be wedge product with $\omega_0$.
\begin{remark}\label{remark multiplic with omega_0 is a morphism of complexes}
Because $\bar\6 \omega_0 = 0$, for $p\geq 0$ we have that $L_{\omega_0}$ is a morphism of complexes 
\[ L_{\omega_0}: (\Lambda^{p,*}, \bar\6) \rightarrow (\Lambda^{p+1, *+1}, \bar\6)  \]
\end{remark}

The restriction and corestriction of $L_{\omega_0}$ to invariant forms,
\[ {L_{\omega_0}: \Lambda^*(M)^{\inv} \rightarrow \Lambda^*(M)^{\inv}} \]
is well defined because $\Lie_{\xi_1} \omega_0 = \Lie_{\xi_2} \omega_0 = 0$. 
In fact, 
$L_{\omega_0}$ is a well defined morphism $L_{\omega_0}: \Lambda_{B, \eta^{0,1}}^{p,q} \rightarrow \Lambda_{B, \eta^{0,1}}^{p+1, q+1}$, which follows because $\omega_0$ is a basic $(1,1)$-form, so whenever $\beta \in (\Lambda^{p,q})_{\bas}$, then $\omega_0 \wedge \beta \in (\Lambda^{p+1,q+1})_{\bas}$. 
Together with Remark \ref{remark multiplic with omega_0 is a morphism of complexes}, this shows that for each fixed $p \geq 1$,
$L_{\omega_0}$ is a morphism of complexes 
\[ L_{\omega_0}: (\Lambda_{B, \eta^{0,1}}^{p-1,*}, \bar\6) \longrightarrow (\Lambda_{B, \eta^{0,1}}^{p, *+1}, \bar\6) \]

Note also that $\bar\6$ takes invariant form to invariant forms since if $\beta$ is an invariant $(p,q)$-form then $0 = d \Lie_{\xi_i} \beta = \Lie_{\xi_i} d \beta = \Lie_{\xi_i} \6 \beta + \Lie_{\xi_i} \bar\6 \beta$ and so $\Lie_{\xi_i} \bar\6 \beta = 0$ because $\Lie_{\xi_i} \bar\6 \beta \in \Lambda^{p,q+1}$ and $\Lie_{\xi_i} \6\beta \in \Lambda^{p+1,q}$. 

Recall the following definition:
\begin{definition}
    Let $(C^*, d_C), (D^*, d_D)$ be complexes and $f: C^* \rightarrow D^*$ be a morphism of complexes. The \textbf{cone of the morphism $f$} is defined to be the complex $(C(f), d_f)$ with $C(f)_i:= C_{i+1} \oplus D_i$ and for $c \in C_{i+1}, d \in D_i$, $d_f( c, d):= (d_C(c), f(c) - d_D(d))$.
\end{definition}

\begin{lemma}\label{lemma complex of invariant forms iso to cone}
    For each fixed $p \geq 0$, the complex $((\Lambda^{p,*}(M))^\inv, \bar \6)$ is isomorphic to the cone of 
    \[ L_{\omega_0}: (\Lambda_{B, \eta^{0,1}}^{p-1,*}, \bar\6) \longrightarrow (\Lambda_{B, \eta^{0,1}}^{p, *+1}, \bar\6) \] shifted by $-1$ i.e. to $C(L_{\omega_0})[-1]$.
    \begin{proof}
        Forms on the tangent space of the foliation $\langle \xi_1, \xi_2 \rangle$ are spanned by $\eta_1, \eta_2$, and hence by $\eta^{0,1}, \eta^{1,0}$. Therefore 
        \begin{align*}
        (\Lambda^{p,q})^\inv &=  (\Lambda^{p,q}_\bas) \oplus
        \left( \Lambda^{p-1,q}_\bas \wedge \eta^{1,0} \right)
        \oplus
        \left( \Lambda^{p,q-1}_\bas \wedge \eta^{0,1} \right) 
        \oplus
        \left(  \Lambda^{p-1,q-1}_\bas \wedge \eta^{0,1} \wedge \eta^{1,0} \right) \\
        &= \Lambda^{p,q}_{B, \eta^{0,1}} \oplus 
        \left( \Lambda^{p-1,q}_{B, \eta^{0,1}} \wedge \eta^{1,0} \right)
        \end{align*}

        The differential $\bar \6$ acts on $\left( \Lambda^{p-1,q}_{B, \eta^{0,1}} \wedge \eta^{1,0} \right)$ as $\bar\6_\bas + L_{\omega_0}$, where
        \[
        {\bar\6_\bas: \Lambda^{p-1,q}_{B, \eta^{0,1}} \wedge \eta^{1,0} \rightarrow \Lambda^{p-1,q+1}_{B, \eta^{0,1}} \wedge \eta^{1,0}}
        \]
        is $\bar\6$ applied to the $\Lambda^{p-1,q}_{B, \eta^{0,1}}$ part, while $L_{\omega_0}$ is multiplication of forms in $\Lambda^{p-1,q}_{B, \eta^{0,1}}$ with $\bar\6 \eta^{1,0} = \omega_0$. 

        \begin{sloppypar}
        This suggests seeing the complex $(\Lambda_{B, \eta^{0,1}}^{p-1,*}, \bar\6)$ as identified with ${(\Lambda_{B, \eta^{0,1}}^{p-1,*}\wedge \eta^{1,0}, \bar\6_\bas)}$; this identification is immediately obtained by simply dropping $\eta^{1,0}$. Seeing $L_{\omega_0}$ after this identification as a morphism of complexes
        \[ L_{\omega_0}: {(\Lambda_{B, \eta^{0,1}}^{p-1,*}\wedge \eta^{1,0}, \bar\6_\bas)} \longrightarrow (\Lambda_{B, \eta^{0,1}}^{p, *+1}, \bar\6) , \]
         the cone of $L_{\omega_0}$ is in degree $q-1$:
        \end{sloppypar}
        \begin{align*}
            \left( C(L_{\omega_0})[-1] \right)_q = \left( C(L_{\omega_0}) \right)_{q-1} 
            = \left( \Lambda^{p-1,q}_{B, \eta^{0,1}} \wedge \eta^{1,0} \right) \oplus \Lambda^{p,q}_{B, \eta^{0,1}}
        \end{align*}

        Thus $\left( C(L_{\omega_0})[-1] \right)_q = \left( \Lambda^{p,q} \right)^\inv$. 

        At position $q-1$ of the cone, the cone differential takes an $\alpha \wedge \eta^{1,0} \in \left( \Lambda^{p-1,q}_{B, \eta^{0,1}} \wedge \eta^{1,0} \right)$ and a $\beta \in \Lambda^{p,q}_{B, \eta^{0,1}}$ to 
        \[ \left( \bar\6_{\bas} (\alpha \wedge \eta^{1,0}),\ L_{\omega_0}(\alpha \wedge \eta^{1,0}) - \bar\6 \beta \right) \]
        Now
        \[ \bar\6_{\bas} (\alpha \wedge \eta^{1,0}) = \left( \bar\6 \alpha \right) \wedge \eta^{1,0} \]
        and by the identification, ${L_{\omega_0}(\alpha \wedge \eta^{1,0}) = \omega_0 \wedge \alpha \in \Lambda^{p,q+1}_{B, \eta^{0,1}}}$. Therefore, the action of the differential of the cone is precisely the same as that of $\bar\6$ and the complex of invariant forms is identified with the $-1$ shift of the cone of $L_{\omega_0}$.
    \end{proof}
\end{lemma}

Furthermore, whenever a compact group acts by holomorphic isometries on a Hermitian manifold, its action on Dolbeault cohomology is trivial:

\begin{theorem}{ \upshape \cite[Theorem~3.3]{klemyatin} }\label{theorem trivial action of group on dolbeault cohom}
    Let $G$ be a compact Lie group acting on a compact Hermitian manifold $M$ by holomorphic isometries. Then the action of $G$ on Dolbeault cohomology, given by $g \cdot [\alpha]:= [g^* \alpha]$ for $g \in G$ and $[\alpha] \in H^{p,q}_{\bar\6}(M)$, is trivial.
\end{theorem}

Consider the unique bi-invariant top form $\nu$ on the compact Lie group $K$ (defined above) with $\int_K \nu = 1$. For any $\alpha \in \Lambda^*(M)$ consider $\overline \alpha := \int_K (k^* \alpha) d\nu(k)$. Then $\bar\alpha$ is an invariant (\cite[Proposition~13.11]{tu_intro_equiv}) smooth (\cite[Proposition~13.13]{tu_intro_equiv}) form of the same degree as $\alpha$. By Theorem \ref{theorem trivial action of group on dolbeault cohom}, taking $\alpha$ to be $\bar\6$-closed, we have for some forms $\beta(k)$
\begin{align*}
    \int_K (k^* \alpha) d\nu(k) &= \int_K \left( \alpha + \bar\6 \beta(k)  \right) d\nu(k) \\
    &= \alpha + \int_K ( \bar\6 \beta(k) ) d\nu(k) = \alpha + \bar\6 \left( \int_K \beta(k) d\nu(k) \right)
\end{align*}
Hence, the cohomology groups $H_{\bar\6}^{p,q}(M)$ are the same as the cohomology groups of $(\Lambda^{p,*}(M)^\inv, \bar\6)$, and hence, by Lemma \ref{lemma complex of invariant forms iso to cone},  
\begin{equation}\label{equation dolbeault cohom same as cohom of cone}
H_{\bar\6}^{p,q}(M) = H^q \left( \left(C \left(L_{\omega_0}:\Lambda_{B, \eta^{0,1}}^{p-1,*} \rightarrow \Lambda_{B, \eta^{0,1}}^{p, *+1} \right) \right) [-1] \right). 
\end{equation}
Now we can prove the following theorem (which has an analogue in the Vaisman setting, \cite[Theorem~4.12]{klemyatin}).

\begin{theorem}
    Let $M$ be the product of two compact Sasakian manifolds with complex structure given by \eqref{explicit formula for pi}. The Dolbeault cohomology groups of $M$ are computed as:
    \begin{align*}
        H_{\bar\6}^{p,q}(M) = \begin{cases}
            \displaystyle \frac{H_{\bas}^{p,q} \oplus [\eta^{0,1}] \wedge H_{\bas}^{p,q-1}(M)}{\im (L_{\omega_0})}, \qquad p+q \leq \dim_{\C} (M) \\[1em]
            \ker (L_{\omega_0}) \rvert_{H_{\bas}^{p,q} \oplus [\eta^{0,1}] \wedge H_{\bas}^{p,q-1}(M)}, \qquad p+q > \dim_{\C} (M)
        \end{cases}
    \end{align*}

    \begin{proof}
        The cone of the morphism $L_{\omega_0}$ gives a short exact sequences of complexes:
        \[ 
            0 \longrightarrow 
             \Lambda_{B, \eta^{0,1}}^{p, *+1}
            \longrightarrow 
            C(L_{\omega_0})
            \longrightarrow
            \left( \Lambda_{B, \eta^{0,1}}^{p-1, *} \right)[1]
        \]
        which gives rise to a long exact sequence in cohomology with connecting map $L_{\omega_0}$:
        \begin{align*}
        \cdots \longrightarrow
        H^{i-1}_{\bar\6} \left( \left( \Lambda_{B, \eta^{0,1}}^{p-1, *} \right)[1] \right)
        &\xrightarrow{L_{\omega_0}}
        H^{i}_{\bar\6} \left( \Lambda_{B, \eta^{0,1}}^{p, *+1} \right)
        \longrightarrow
        H^{i} ( C(L_{\omega_0}) )
        \longrightarrow \\
        & \longrightarrow
        H^{i}_{\bar\6} \left( \left( \Lambda_{B, \eta^{0,1}}^{p-1, *} \right)[1] \right)
        \xrightarrow{L_{\omega_0}} \cdots
        \end{align*}

        Taking into account shifts, degrees and \eqref{equation dolbeault cohom same as cohom of cone} we thus have:
         \begin{align*}
        \cdots \longrightarrow
        H^{i}_{\bar\6} \left( \Lambda_{B, \eta^{0,1}}^{p-1, *}\right)
        &\xrightarrow{L_{\omega_0}}
        H^{i}_{\bar\6} \left( \Lambda_{B, \eta^{0,1}}^{p, *+1} \right)
        \longrightarrow
        H^{p, i+1}_{\bar\6}(M)
        \longrightarrow  \\
        &\longrightarrow
        H^{i+1}_{\bar\6} \left(  \Lambda_{B, \eta^{0,1}}^{p-1, *}  \right)
        \xrightarrow{L_{\omega_0}} \cdots
        \end{align*}

        Now since $\bar\6 \eta^{0,1}=0$, 
        \[ H^i_{\bar \6} (\Lambda^{p-1,*}_{B, \eta^{0,1}}) = H^{p-1,i}_{\bas}(M) \oplus [\eta^{0,1}] \wedge H^{p-1, i-1}_{\bas}(M) \]
        By Theorem \ref{theorem trasnversally kaehler admits hodge}, basic cohomology behaves just like the cohomology of a K\"ahler manifold with K\"ahler form $\omega_0$. Hence, since $M$ is compact, by the Hodge isomorphism theorem and the fact that the K\"ahler form is harmonic, the operator
        \[H^{t}_{\bar\6} \left( \Lambda_{B, \eta^{0,1}}^{s, *}\right)
        \xrightarrow{L_{\omega_0}}
        H^{t}_{\bar\6} \left( \Lambda_{B, \eta^{0,1}}^{s+1, *+1} \right) \]
        is injective whenever $s+t \leq \dim_{\C} M - 1$; by Poincar\'e duality, it is surjective whenever $s+t > \dim_{\C} M - 1$.
        Hence, for $p+i \leq \dim_\C M$, we obtain the short exact sequence
        \begin{equation}\label{eqn_short_exact_injective}
        0 \longrightarrow 
        H^{i-1}_{\bar\6} \left( \Lambda_{B, \eta^{0,1}}^{p-1, *}\right)
        \xrightarrow{L_{\omega_0}}
        H^{i}_{\bar\6} \left( \Lambda_{B, \eta^{0,1}}^{p, *} \right)
        \longrightarrow
        H^{p, i}_{\bar\6}(M)
        \longrightarrow 0
        \end{equation}
        while for $p+i > \dim_\C M + 1$ we obtain the short exact sequence
        \begin{equation}\label{eqn_short_exact_surjective}
        0 \longrightarrow
        H^{p, i}_{\bar\6}(M)
        \longrightarrow 
        H^{i}_{\bar\6} \left(  \Lambda_{B, \eta^{0,1}}^{p-1, *}  \right)
        \xrightarrow{L_{\omega_0}}
        H^{i+1}_{\bar\6} \left(  \Lambda_{B, \eta^{0,1}}^{p, *}  \right)
        \longrightarrow 0
        \end{equation}

        \begin{sloppypar}Finally, when $p+i = \dim_\C M +1$, by the Hard Lefschetz theorem
        $\displaystyle {H^{i-1}_{\bar\6} \left( \Lambda_{B, \eta^{0,1}}^{p-1, *}\right)
        \xrightarrow{L_{\omega_0}}
        H^{i}_{\bar\6} \left( \Lambda_{B, \eta^{0,1}}^{p, *} \right)}$
        is an isomorphism, so in particular surjective.
        Since, as mentioned above, $\displaystyle {H^{i}_{\bar\6} \left( \Lambda_{B, \eta^{0,1}}^{p-1, *}\right)
        \xrightarrow{L_{\omega_0}}
        H^{i+1}_{\bar\6} \left( \Lambda_{B, \eta^{0,1}}^{p, *} \right)}$
        is also surjective, we have the short exact sequence \eqref{eqn_short_exact_surjective} also for the case when $p+i = \dim_\C M+1$.
        \end{sloppypar}
    \end{proof}
\end{theorem}

\hfill

{\small
	
	\noindent {\sc Vlad Marchidanu\\
		University of Bucharest, Faculty of Mathematics and Informatics, \\14
		Academiei str., 70109 Bucharest, Romania\\
		\tt marchidanuvlad@gmail.com}
}
\end{document}